\RequirePackage{fix-cm}

\documentclass[smallextended]{my_svjour3}

\smartqed

\usepackage{amsmath}
\usepackage{graphicx}
\usepackage[update,prepend]{epstopdf}
\usepackage{url}

\begin{document}

\journalname{Math. Comput. Sci.}

\date{Submitted 16-July-2015 / Accepted, after a revision, 10-Mar-2016.} 

\title{Dynamics and optimal control of Ebola transmission}

\titlerunning{Dynamics and optimal control of Ebola transmission}

% ------------------------------------------------------

\authorrunning{A. Rachah, D. F. M. Torres}

\author{Amira Rachah \and Delfim F. M. Torres}

% ------------------------------------------------------

\institute{Amira Rachah \at
Math\'{e}matiques pour l'Industrie et la Physique,\\
Institut de Math\'{e}matiques de Toulouse, Universit\'e Paul Sabatier,\\
F-31062 Toulouse Cedex 9, France\\
\email{arachah@math.univ-toulouse.fr}
\and
Delfim F. M. Torres (corresponding author) \at
Center for Research and Development in Mathematics and Applications (CIDMA), \\
Department of Mathematics, University of Aveiro, 3810-193 Aveiro, Portugal \\
\email{delfim@ua.pt}
}

\maketitle

% ------------------------------------------------------

\begin{abstract}
A major Ebola outbreak occurs in West Africa since March 2014, 
being the deadliest epidemic in history. As an infectious 
disease epidemiology, Ebola is the most lethal 
and is moving faster than in previous outbreaks.
On 8 August 2014, the World Health Organization (WHO) declared
the outbreak a public health emergency of international concern.
Last update on 7 July 2015 by WHO reports 27 609 cases of Ebola
with a total of 11 261 deaths. In this work, we present 
a mathematical description of the spread of Ebola virus based 
on the SEIR (Susceptible--Exposed--Infective--Recovered) model 
and optimal strategies for Ebola control. In order to control 
the propagation of the virus and to predict the impact of vaccine programmes,
we investigate several strategies of optimal control of the spread of Ebola: 
control infection by vaccination of susceptible; 
minimize exposed and infected; reduce Ebola
infection by vaccination and education.

\keywords{Mathematical modelling \and optimal control \and epidemiology 
\and Ebola virus \and 2014 Ebola epidemic in West Africa.}
\end{abstract}

% ------------------------------------

\section{Introduction}

The most widespread epidemic of Ebola virus disease in history is currently 
ongoing in West Africa. Ebola virus is the most complex and lethal pathogen 
for humans \cite{barraya,joseph}. The first outbreak was in $1976$ in Congo, 
close to Ebola river, where the disease takes its name. The Ebola virus causes 
an acute, often fatal haemorrhagic illness. The current outbreak in West Africa 
is the largest and most complex Ebola outbreak since the virus discovery.
The increasing of this virulent virus is extremely rapid and has evolved 
into a health and humanitarian catastrophe of historic scope. The incubation 
period, that is, the time interval from infection with the virus to the 
onset of symptoms, is $2$ to $21$ days. Humans are not infectious to another 
person until they develop symptoms. Symptoms typically begin with sudden onset 
of  fever, sore throat, muscle pain, and headaches. Then, vomiting, diarrhoea 
and rash usually follow, along with decreased function of the liver and kidneys.
At this time infectious individuals begin to bleed, both internally and externally 
(bleeding from nose, mouth, eyes and anus). In final stage, these symptoms are
followed by death \cite{borio,dowel,legrand,peter,anon2,okwar,anon1}. Ebola 
spreads through human-to-human transmission via close and direct physical contact
(through broken skin or mucous membranes) with infected bodily fluids.
The most infectious fluids are blood, faeces and vomit secretions; however, 
all body fluids have the capacity to transmit the virus. The virus is also 
transmitted indirectly via exposure to objects or environment contaminated 
with infected secretions. Because of this, healthcare workers must practice 
strict infection prevention and control precautions.

In epidemiology, mathematical models are a key tool that contribute 
to the understanding of the dynamics of a virus and the impact of 
vaccination programmes. More precisely, Mathematics has an important role 
in the study of propagation of virus spreads by allowing policy-makers 
to predict the impact of particular vaccine programmes or to derive more 
efficient strategies based on mathematical insights \cite{delf,MR3101449}. 
In recent years, optimal control theory has become a powerful mathematical tool 
to assess the intervention of public health authorities. Indeed, the inclusion, 
in an epidemic model, of some practical control strategies, like vaccines, 
social distancing or quarantine, provides a rational basis for policies designed 
to control the spread of the virus. In this spirit, the aim of this work is to  
investigate effective strategies to control the spread of the Ebola virus 
by setting appropriate optimal control problems subject to a SEIR 
epidemic model that divides the population into four groups: the Susceptible,
the Exposed (infected but not infectious), the Infectious, and the Recovered 
\cite{althaus,astacio,Gerard,chow,zeng}. For a comparison study between 
the SIR and SEIR models we refer the reader to \cite{MyID:331}.
The main practical strategy of optimal control adopted is a vaccine
\cite{vacc_opt2,vacc_opt4,vacc_opt3,vacc_opt5,vacc_opt1}. 
For fractional (non-integer order) models see \cite{Area1,Area2}.

The paper is organized as follows. In Section~\ref{sec:2} we present
the mathematical model that describes the dynamics of the spread of Ebola virus
in West Africa. In Section~\ref{sec:3} we study several control strategies  
for the propagation of the virus by using the proposed model. In these strategies, 
we use parameters estimated from recent statistical data based
on the WHO report of the 2014 Ebola outbreak \cite{who}. 
In Sections~\ref{subsec:3.1} and \ref{subsec:3.2} we study the case 
of control through vaccination treatment by considering two different
objective functions. In Section~\ref{subsec:3.3} we consider educational 
campaigns as a control, coupled with the vaccination treatment. 
We end with Section~\ref{subsec:3.5} of discussion of the numerical results 
and with Section~\ref{sec:4} of conclusions and future work.

% ------------------------------------

\section{The basic model equations}
\label{sec:2}

\begin{figure}
\centering
\includegraphics[width=10cm]{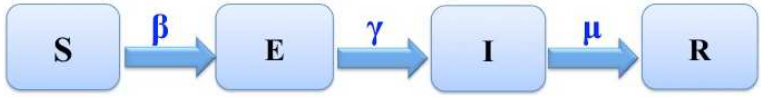}
\caption{Compartment diagram of the 
Susceptible--Exposed--Infectious--Recovered (SEIR) model \eqref{eq1:SEIR}.}
\label{SEIR_fig1}
\end{figure}

In this section, we describe the transmission of Ebola virus by a SEIR model,  
where the population is divided into four groups: the susceptible individuals 
at time $t$, denoted by $S(t)$, enter the exposed class $E(t)$ before
they become infectious. The infectious class at time $t$, denoted by $I(t)$,
represents the individuals that are infected with the disease and are suffering 
the symptoms of Ebola. Finally, we have the recovered class, which at time $t$ 
is denoted by $R(t)$. The total population, assumed constant during the short 
period of time under study, is given by 
\begin{equation}
\label{eq:N}
N = S(t) + E(t) + I(t) + R(t)
\end{equation}
at any instant of time $t$. The transmission of Ebola virus is then described 
by the following set of nonlinear ordinary differential equations (ODEs):
\begin{equation}
\label{eq1:SEIR}
\begin{cases}
\dfrac{dS(t)}{dt} = -\beta S(t)I(t),\\[0.30cm]
\dfrac{dE(t)}{dt} = \beta S(t)I(t) - \gamma E(t),\\[0.30cm]
\dfrac{dI(t)}{dt} = \gamma E(t) - \mu I(t),\\[0.30cm]
\dfrac{dR(t)}{dt} = \mu I(t).
\end{cases}
\end{equation}
Transitions between different states are described by the following parameters:
\begin{itemize}
\item the transmission rate $\beta$,
\item the infectious rate $\gamma$,
\item the recovered rate $\mu$.
\end{itemize}
Figure~\ref{SEIR_fig1} shows the relationship between 
the variables of model \eqref{eq1:SEIR}. Its numerical 
resolution is studied by Rachah and Torres
in \cite{symcomp}. See also \cite{MyID:321} 
for a simpler SIR model describing Ebola.

% ----------------------------------------

\begin{remark}
The model we consider here is very simple and assumes that the population 
is constant in the period of time under study: the sum of the right-hand 
side of the equations of system \eqref{eq1:SEIR} is zero. The
possibility to design more complex models for a non-constant
population is possible: see the recent eight state-variables model 
of \cite{Area2}. However, for the analysis and simulation of the model 
to be possible, the authors of \cite{Area2} imposed conditions 
to maintain the population $N$ constant. As remarked 
at the end of \cite{Area2}, the design of strategies 
when the population is not constant seems 
to be a good topic of research. 
\end{remark}

% ----------------------------------------

\section{Strategies for the control of the virus}
\label{sec:3}

Recently, epidemiological models have used optimal control techniques. 
Most works focus on HIV, tuberculosis (TB) and dengue  
\cite{vacc_opt4,vacc_opt3,MR2719552,delf,MR3266821,MR3101449}.
Control efforts are carried out to limit the spread of the disease 
and, in some cases, to prevent the emergence of drug resistance. In this section, 
we formulate several optimal control problems subject to the SEIR model 
\eqref{eq1:SEIR}, in order to derive the optimal treatment strategies. 
For each strategy, we study a specific objective in order to minimize not only 
the number of infected individuals and the systemic costs of treatments, but 
also to include an educational campaign with the vaccination treatment. 
The integration of educational campaigns has a great importance in
countries that do not have the capacity to defend themselves against the virus.
We compare the result of each strategy with the simulation results
previously studied by Rachah and Torres in \cite{symcomp}, the so called Strategy~1,
which is described in Section~\ref{subsec:3.1}. Strategies~2
and $3$ are an improvement of the Strategy~1 of \cite{symcomp}, 
and are given in Sections~\ref{subsec:3.2} and \ref{subsec:3.3}, respectively. 
More precisely, Strategy~3 consists in the study of an educational campaign 
(about the virus) coupled with a vaccination treatment. A comparison between 
the different strategies and the simulation results 
is given in Section~\ref{subsec:3.5}.

Two approaches are common in optimal control: 
direct and indirect methods \cite{MyID:287}. 
The direct methods  are based in the discretization 
of the optimal control problems, reducing them 
to nonlinear constrained optimization problems. 
Indirect methods are based on the Hamiltonian
and the Pontryagin maximum principle. 
In this work we have adopted the direct approach.
We have computed the numerical solutions of the optimal control problems 
with the help of the \textsf{ACADO} solver \cite{acado}. \textsf{ACADO} 
is based on a multiple shooting method, including automatic differentiation 
and based ultimately on the semidirect multiple shooting algorithm 
of Bock and Pitt \cite{acado2}. The \textsf{ACADO} solver comes 
as a self-contained public domain software environment, written in \textsf{C++}, 
for automatic control and dynamic optimization.
The reader interested in the details of the implemented algorithms 
is referred to \cite{MR2839653}.

% ------------------------------------

\subsection{Strategy 1}
\label{subsec:3.1}

In this subsection, we present the optimal control problem 
investigated in \cite{symcomp}, obtained by introducing 
into the model \eqref{eq1:SEIR} a control $u(t)$ representing the vaccination 
rate at time $t$. The control $u(t)$ is the fraction of susceptible individuals 
being vaccinated per unit of time. Then, the mathematical model with control 
is given by the following system of nonlinear differential equations:
\begin{equation}
\label{SEIR_control}
\begin{cases}
\dfrac{dS(t)}{dt} = -\beta S(t)I(t) - u(t) S(t),\\[0.30cm]
\dfrac{dE(t)}{dt} = \beta S(t)I(t) - \gamma E(t),\\[0.30cm]
\dfrac{dI(t)}{dt} = \gamma E(t) - \mu I(t),\\[0.30cm]
\dfrac{dR(t)}{dt} = \mu I(t) + u(t) S(t).
\end{cases}
\end{equation}
The goal of the strategy is to reduce the infected individuals and the cost
of vaccination. Precisely, the optimal control problem consists of minimizing
the objective functional
\begin{equation}
\label{cost_func_strat1}
J(u) = \int_{0}^{t_{end}} \left[I(t) + \dfrac{\tau}{2}u^2(t)\right] dt,
\end{equation}
where $u(t)$ is the control variable, which represents the vaccination rate
at time $t$, and the positive parameters $\tau$ and $t_{end}$ denote, respectively,
the weight on cost and the duration of the vaccination program.

% -----------------------------------

\subsection{Strategy 2}
\label{subsec:3.2}

We now improve the strategy of Section~\ref{subsec:3.1}   
by studying other strategy in order to better control 
the propagation of the spread of Ebola into populations.
Our goal in this strategy is to reduce the number of exposed
and infected individuals. More precisely, our optimal control 
problem consists of minimizing the objective functional
\begin{equation}
\label{cost_func_strat2}
J(u) = \int_{0}^{t_{end}}
\left[ A_1E(t) + A_2 I(t) + \dfrac{\nu}{2}u^2(t) \right] dt
\end{equation}
subject to the model described by \eqref{SEIR_control}.
The two first terms in the functional objective \eqref{cost_func_strat2}
represent benefit of $E(t)$ and $I(t)$ populations that we wish to reduce; 
$A_1$ and $A_2$ are positive constants to keep a balance in the size 
of $E(t)$ and $I(t)$, respectively. In the quadratic term of 
\eqref{cost_func_strat2}, $\nu$ is a positive weight parameter 
associated with the control $u(t)$, and the square of the control 
variable reflects the severity of the side effects of the vaccination. 
One has $u \in \mathcal{U}_{ad}$, where
$$
\mathcal{U}_{ad}=\left\{u : u \,  \text{is measurable}, 0
\leq u(t) \leq u_{max}<\infty, \, t\in [0,t_{end}] \right\}
$$
is the admissible control set, with $u_{max}=0.9$.
Note that it is not realistic to admit the possibility
to vaccinate everybody and the value $0.9$ means that, 
at maximum, 90\% of susceptible are vaccinated.
In the numerical simulations of model \eqref{eq1:SEIR} and Strategy~1
of optimal control, Rachah and Torres \cite{symcomp} used
the parameters estimated on November $2014$ by Althaus and Kaurov  
\cite{althaus,kaurov}, who studied statistically the data
of the World Health Organisation (WHO) \cite{who}.
In particular, in the statistical study \cite{kaurov},
Kaurov studies the outbreak by modeling it
with Wolfram's \textsf{Mathematica} language.

\begin{figure}
\centering
\includegraphics[width=10cm]{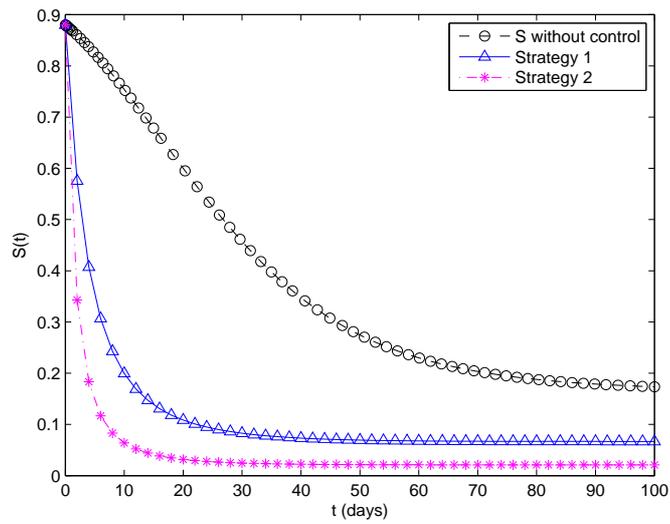}
\caption{Comparison between the curves of susceptible individuals $S(t)$
in case of Strategy~$1$ and Strategy~$2$ \emph{versus} 
without control. \label{stratg12_S}}
\end{figure}

\begin{figure}
\centering
\includegraphics[width=10cm]{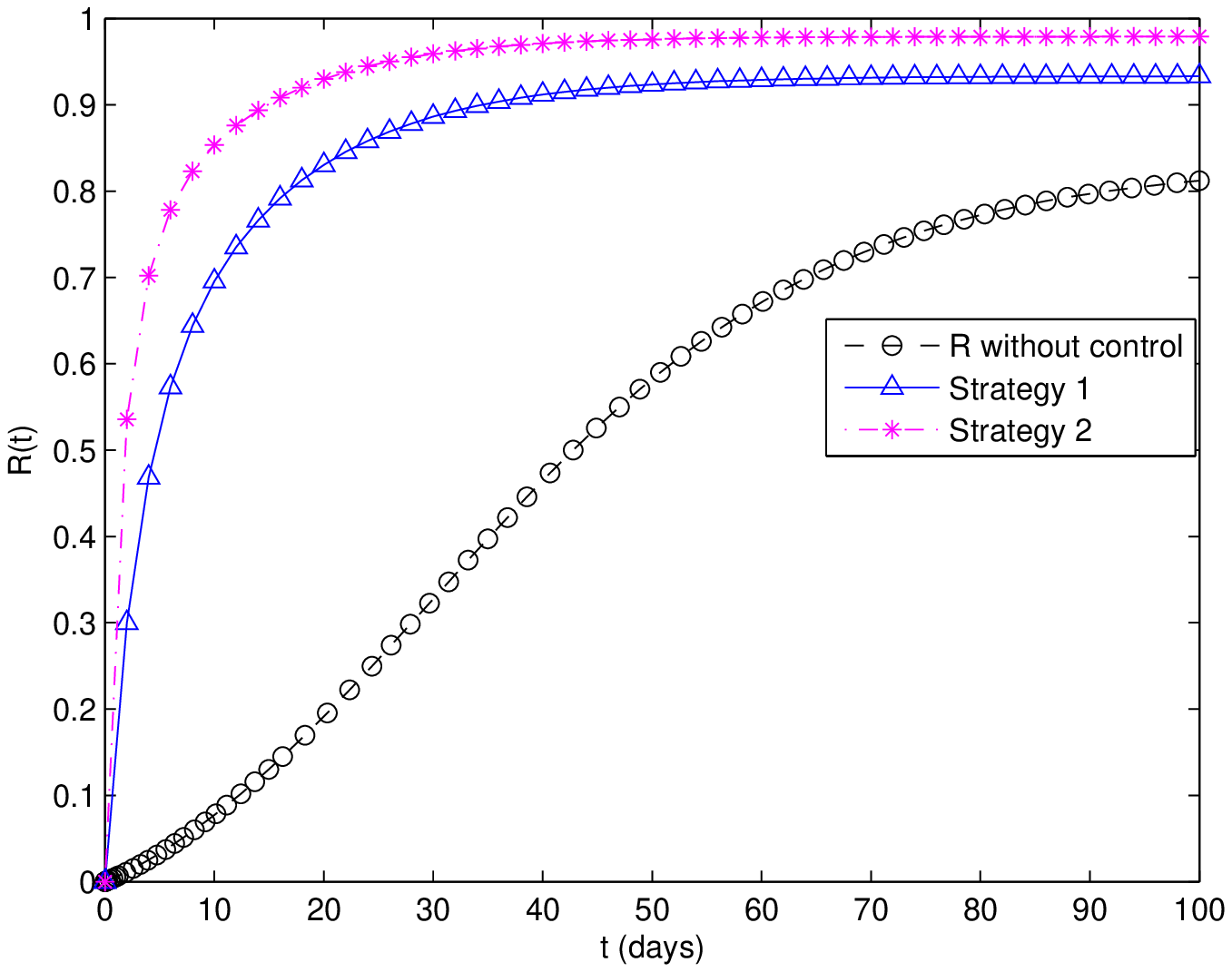}
\caption{Comparison between the curves of recovered individuals $R(t)$
in case of Strategy~$1$ and Strategy~$2$ \emph{versus}
without control. \label{stratg12_R}}
\end{figure}

\begin{figure}
\centering
\includegraphics[width=10cm]{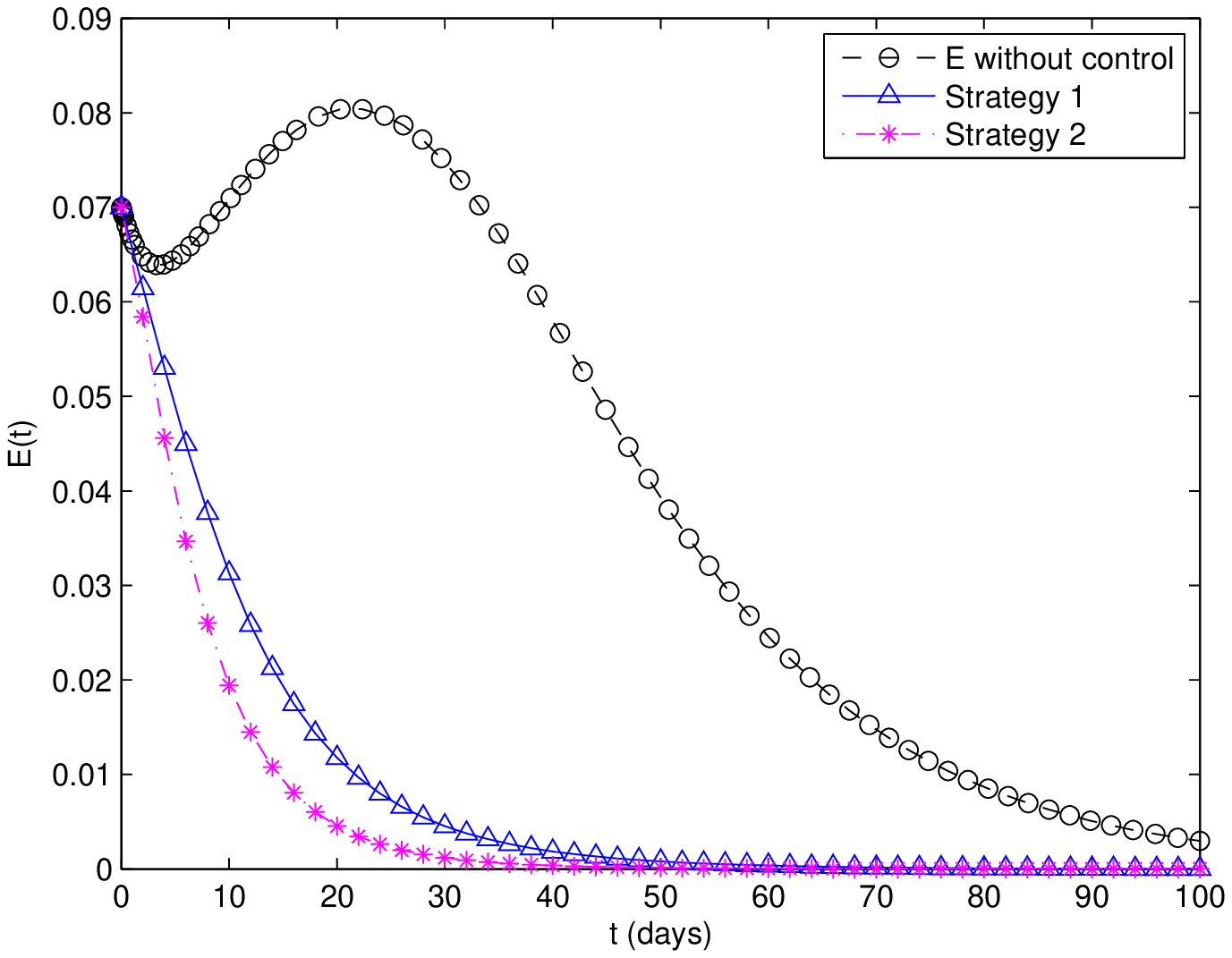}
\caption{Comparison between the curves of exposed individuals $E(t)$
in case of Strategy~$1$ and Strategy~$2$ \emph{versus}
without control. \label{stratg12_E}}
\end{figure}

\begin{figure}
\centering
\includegraphics[width=10cm]{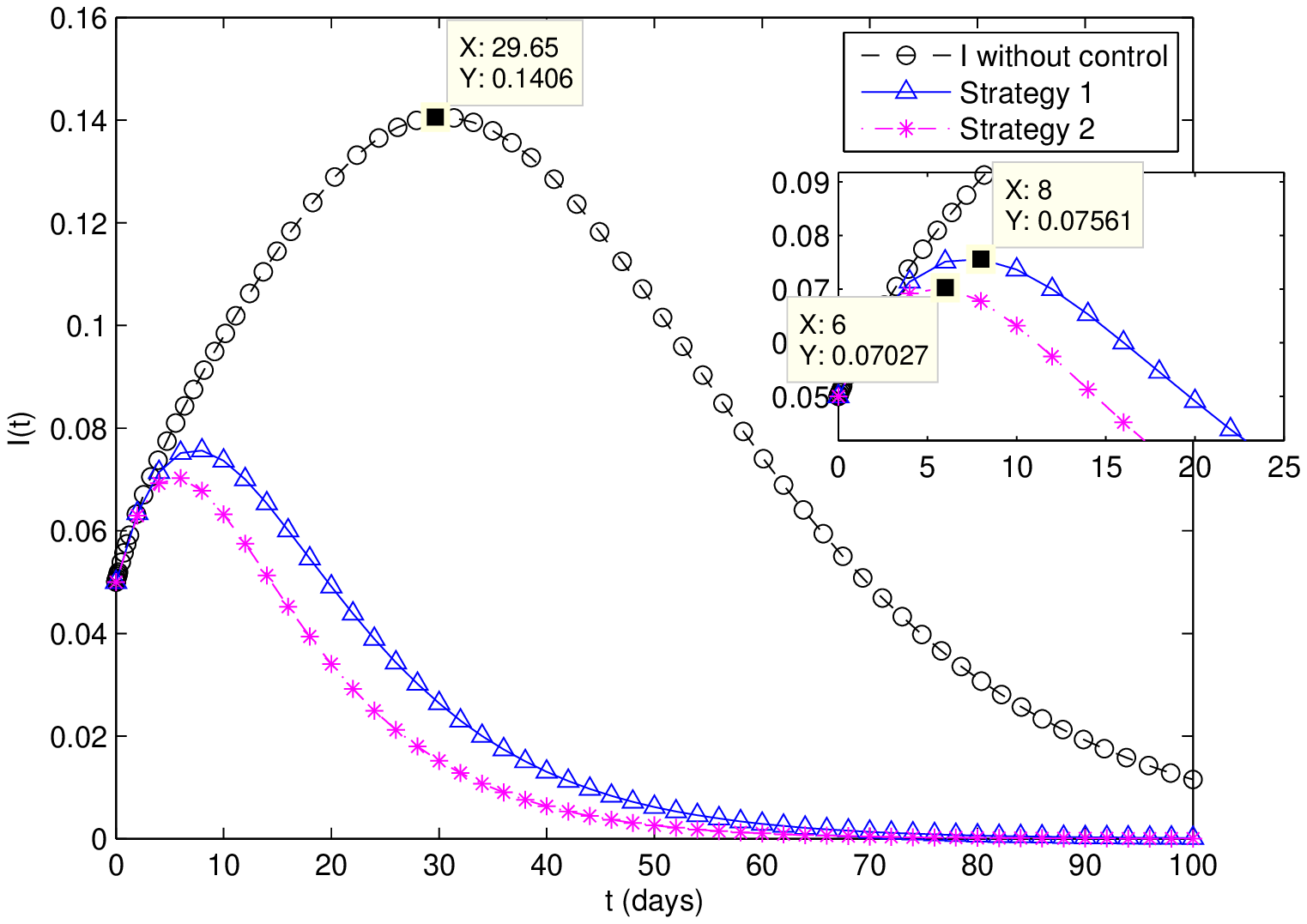}
\caption{Comparison between the curves of infected individuals $I(t)$
in case of Strategy~$1$ and Strategy~$2$ \emph{versus}
without control. \label{stratg12_I}}
\end{figure}

\begin{figure}
\centering
\includegraphics[width=10cm]{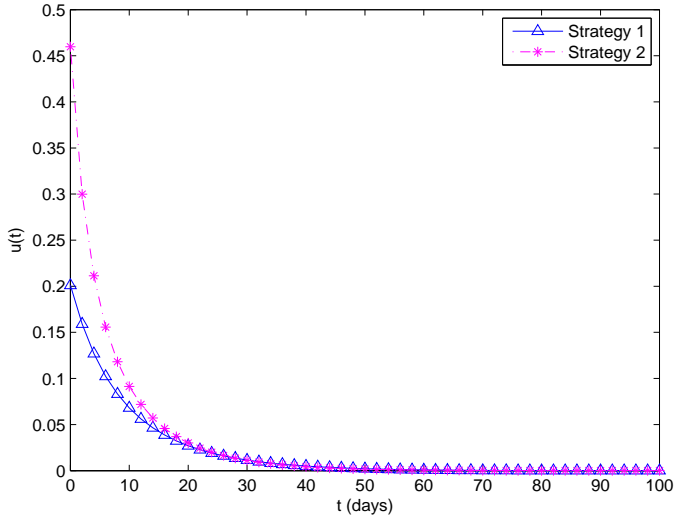}
\caption{The optimal control $u$ for Strategy~$1$ and Strategy~$2$.
\label{stratg2_u}}
\end{figure}

In order to compare our improvement of the optimal control study 
with the previous results of \cite{symcomp}, we use here the same parameters, 
that is, the same rate of infection $\beta=0.2$, the infectious rate 
$\gamma=0.1887$, the same recovered rate $\mu=0.1$, and the same initial values
$\left(S(0),E(0),I(0),R(0)\right) =\left(0.88, 0.07, 0.05, 0\right)$
for the initial number of susceptible, exposed, infected, and recovered populations
(at the beginning, 88\% of population is susceptible, 7\% is exposed
and 5\% is infected with Ebola). Figures~\ref{stratg12_S}, \ref{stratg12_R}, 
\ref{stratg12_E} and \ref{stratg12_I} show, respectively, the significant 
difference in the number of susceptible, recovered, exposed and infected 
individuals with Strategy~$1$, Strategy~$2$, and without control.
In Figure~\ref{stratg12_S}, we see that the number of susceptible $S$,
in case of optimal control under Strategy~$2$, decreases more rapidly during
the vaccination campaign. It reaches 2\% at the end of the campaign,
in contrast with the 6.6\% at the end of the campaign with Strategy~$1$,
and against 17.3\% in the absence of optimal control. Figure~\ref{stratg12_R} 
shows that the number of recovered individuals increases rapidly. The number 
$R(t_{end})$ at the end of the optimal control vaccination period of 
Strategy~2 is 97.9\%, instead of 93.3\% in case of Strategy~$1$, 
and against 81.2\% without control. Figure~\ref{stratg12_E} presents the 
time-dependent curve of exposed individuals. This curve  shows that there is no 
peak of the curve of exposed individuals in case of control with Strategy~1 
and Strategy~2 when we compare it with the curve in case without control, 
which shows a high peak. When we compare Strategy~1 and Strategy~2, we see that 
the curve of exposed individuals decreases more rapidly in case of control 
with Strategy~2. The period of incubation of the virus is only 30 days 
in case of control with Strategy~2 instead of 50 days in  Strategy~1.
In Figure~\ref{stratg12_I}, the time-dependent curve of infected individuals
shows that the peak of the curve of infected individuals is less important
in case of control with Strategy~2. In fact, the maximum value on the infected
curve $I$ under optimal control is 7\% in case of Strategy~2, instead of 7.56\%
in Strategy~$1$, and against 14\% without any control (see Figure~\ref{stratg12_I}).
The other important effect of Strategy~$2$, which we see in the same curve, is 
the period of infection, which is the less important. The value of the period
of infection is $65$ days in case of Strategy~$2$, instead of $80$ days in case 
of Strategy~$1$, and against more than $100$ days without vaccination.
This shows the efficiency of vaccination with Strategy~$2$ in controlling Ebola. 
Figure~\ref{stratg2_u} gives a representation of the optimal control $u(t)$ 
for Strategy~$1$ and Strategy~$2$. The minimum values of the $J(u)$ functionals 
\eqref{cost_func_strat1} and \eqref{cost_func_strat2} are, respectively,
1,936 (Strategy~1) and 2,509 (Strategy~2).

% -------------------------------
%--------------------------------
% -------------------------------
% -------------------------------
% -------------------------------

% -------------------------------

\subsection{Strategy 3}
\label{subsec:3.3}

\begin{figure}
\centering
\includegraphics[width=10cm]{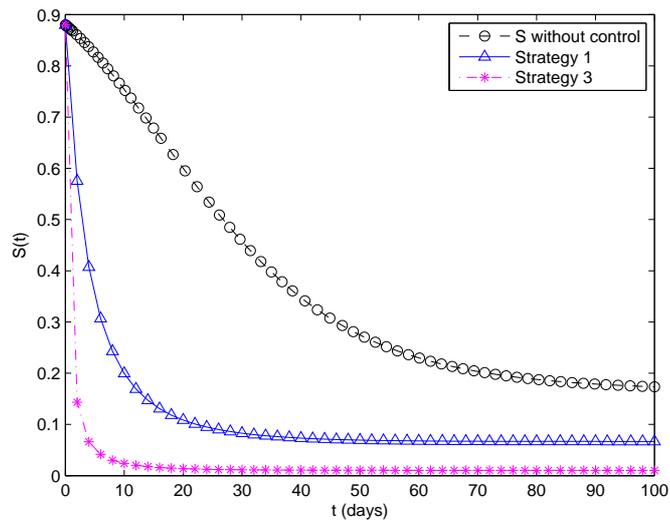}
\caption{Comparison between the curves of susceptible individuals $S(t)$
in case of Strategy~$1$ and Strategy~$3$
\emph{versus} without control. \label{stratg13_S}}
\end{figure}

\begin{figure}
\centering
\includegraphics[width=10cm]{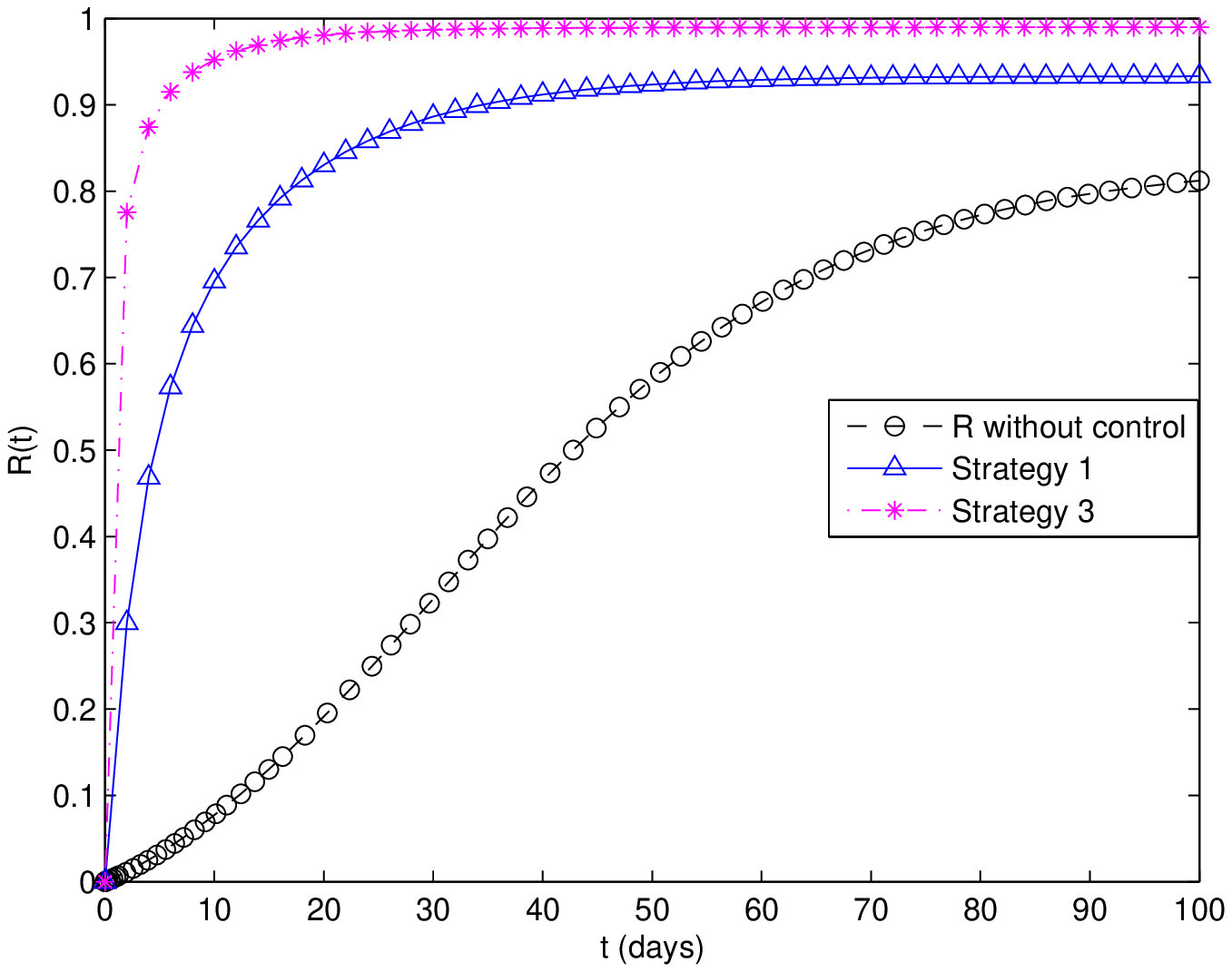}
\caption{Comparison between the curves of recovered individuals $R(t)$
in case of Strategy~$1$ and Strategy~$3$
\emph{versus} without control. \label{stratg13_R}}
\end{figure}

\begin{figure}
\centering
\includegraphics[width=10cm]{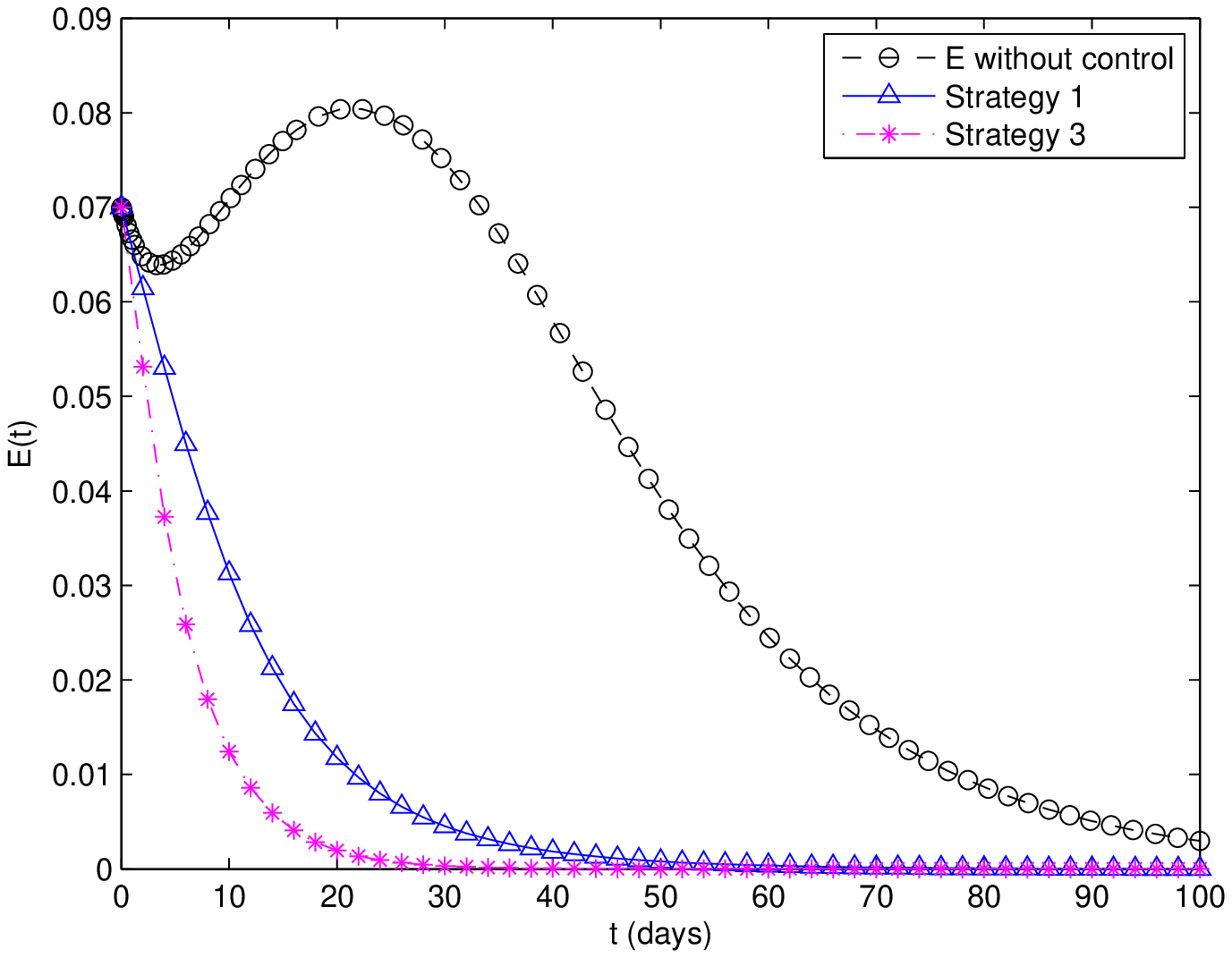}
\caption{Comparison between the curves of exposed individuals $E(t)$
in case of Strategy~$1$ and Strategy~$3$
\emph{versus} without control. \label{stratg13_E}}
\end{figure}

\begin{figure}
\centering
\includegraphics[width=10cm]{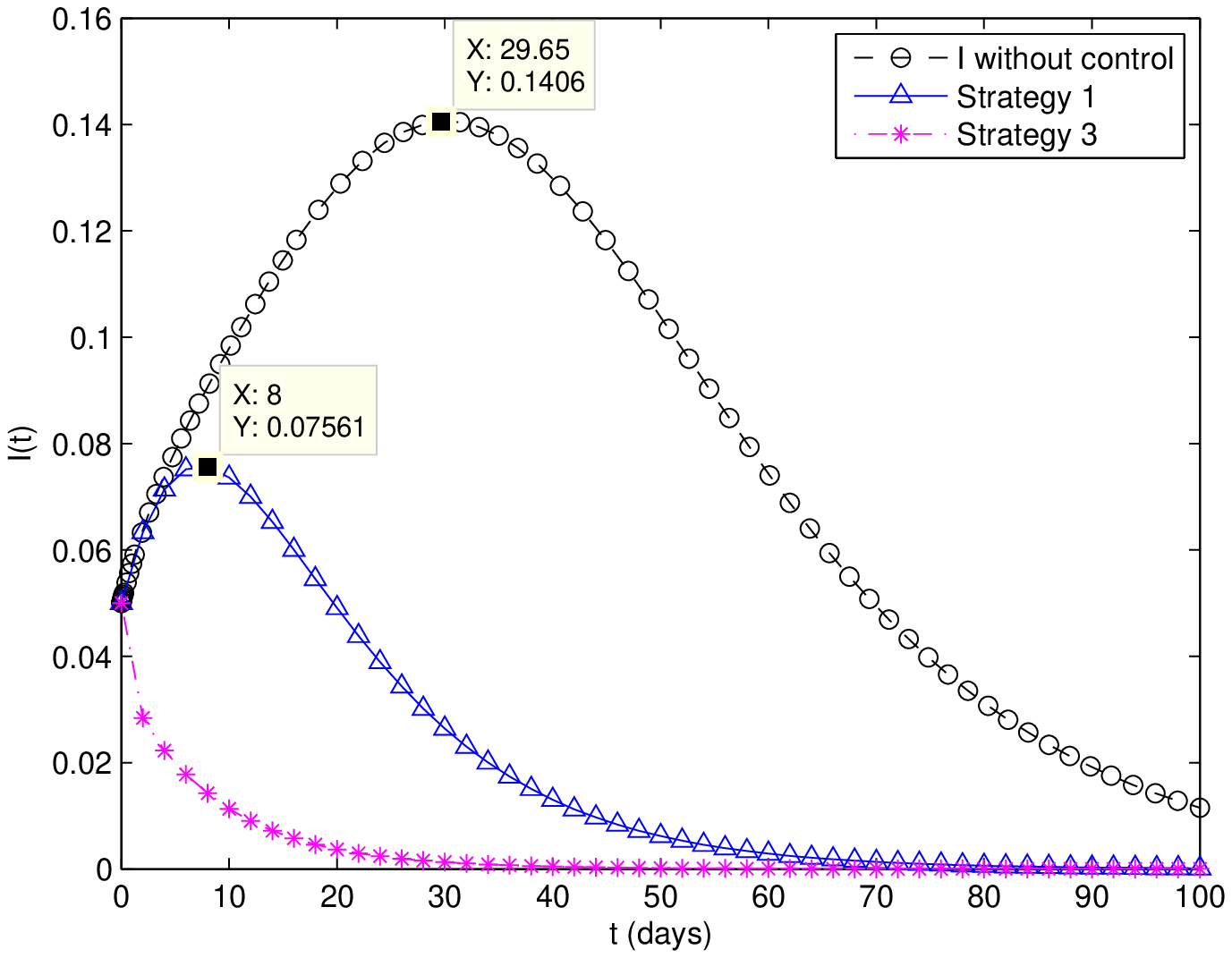}
\caption{Comparison between the curves of infected individuals $I(t)$
in case of Strategy~$1$ and Strategy~$3$
\emph{versus} without control. \label{stratg13_I}}
\end{figure}

\begin{figure}
\centering
\includegraphics[width=10cm]{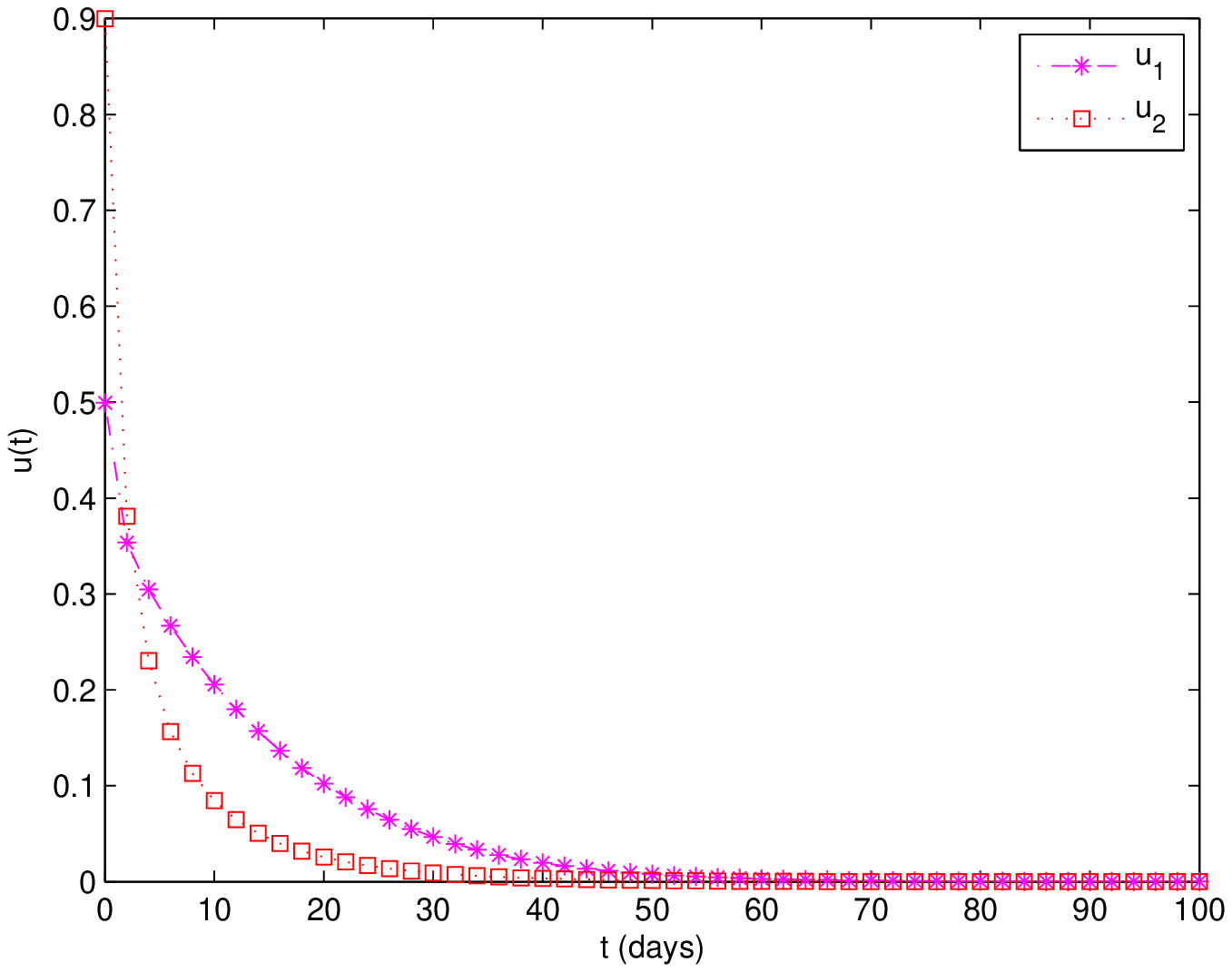}
\caption{The optimal control variables $u_1$ and $u_2$ for Strategy~$3$.
\label{stratg3_u}}
\end{figure}

In this strategy we use the fact that individuals can acquire immunity against
the virus either through educational campaigns or recovery after treatment
for the virus. Our idea is to study the effect of educational campaigns
with a vaccination treatment in practical Ebola situations. The case
of the French nurse cured of Ebola is a proof of the possibility of educational
campaigns and treatment \cite{valler}. An educational campaign, in case
of spread of Ebola, has great importance. In fact, Ebola virus spreads
through human-to-human transmission, not only by close and direct physical
contact with infected bodily fluids, but also via exposure to objects
or contaminated environment. The most infectious fluids are blood,
feces, and vomit secretions. However, all body fluids have the capacity
to transmit the virus. Here, we intend to control the propagation 
of the Ebola virus by using two control variables in the SEIR model, as follows:
\begin{equation}
\label{SEIR_control_stratg3}
\begin{cases}
\dfrac{dS(t)}{dt} = -\beta S(t)I(t) - u_2(t) S(t),\\[0.30cm]
\dfrac{dE(t)}{dt} = \beta S(t)I(t) - \gamma E(t),\\[0.30cm]
\dfrac{dI(t)}{dt} = \gamma E(t) - \mu I(t) - u_1(t)I(t),\\[0.30cm]
\dfrac{dR(t)}{dt} = \mu I(t) + u_1(t) I(t) + u_2(t) S(t),
\end{cases}
\end{equation}
where $u_1(t)$ is the fraction of infective that is treated at time $t$,
and $u_2(t)$  is the fraction of susceptible individuals that is subject
to an educational campaign at time $t$. Our goal is to minimize
simultaneously the total number of individuals that are infected, the cost
of treatment, and the cost of educational campaigns to the population.
The objective functional is now
\begin{equation}
\label{cost_func_strat3}
J(u) = \int_{0}^{t_{end}} \left[ \kappa I(t)
+ B_1\dfrac{u_1^2(t)}{2}+ B_2\dfrac{u_2^2(t)}{2} \right] dt
\end{equation}
subject to system \eqref{SEIR_control_stratg3}, where
$u = (u_1, u_2)$, with $u_1$ representing treatment
and $u_2$ educational campaigns, and $\kappa$, $B_1$ 
and $B_2$ are weight parameters. The Lebesgue
measurable control set is defined as
$$
\mathcal{U}_{ad}:=\left\{u = (u_1,u_2) : u \in L^1,\,
0 \leq u_1(t), \, u_2(t) \leq u_{\max},
\, t\in [0,t_{end}] \right\},
$$
where $u_{max}=0.9$. We choose quadratic terms 
with respect to the controls in order to describe
the nonlinear  behaviour of the cost of implementing the educational campaign
and treatments. The first term in the objective functional \eqref{cost_func_strat3}, 
$\kappa I$, stands for the total number of individuals that are infected; the 
term $B_1 u_1^2 / 2$ represents the cost of treatment; while the term 
$B_2 u_2^2/2$ represents the cost associated with the educational campaign.
Figure~\ref{stratg13_S} shows the time-dependent curve of susceptible
individuals, $S(t)$, which decreases more rapidly in case of the optimal control
with Strategy~$3$. It reaches 1\% at the end of the campaign, instead of 6.6\%
at the end of the campaign in case of Strategy~1, and against 17.3\% in the 
absence of optimal control. Figure~\ref{stratg13_R} shows that the number 
of recovered individuals of Strategy~$3$ increases rapidly until 98.9\%, 
instead of 93.3\% in case of Strategy~$1$, and against 81.2\% without control.
Figure~\ref{stratg13_E} shows that the number of exposed individuals of 
Strategy~$3$ increases the most rapidly. The period of incubation of the virus 
is the lest important (30 days) instead of 60 days in case of Strategy~$1$.
In Figure~\ref{stratg13_I} we see the time-dependent curve of infected
individuals $I(t)$, which decreases mostly. The other important effect
of the Strategy~$3$, which we can see in the same curve, is the period
of infection, which is the less important. The value of the period of 
infection is $30$ days in case of Strategy~$3$, instead of $80$ days 
in case of Strategy~$1$, and against $100$ days without vaccination. 
This shows the efficiency of the effect of educational campaigns 
in controlling Ebola virus with the vaccination control treatment described 
in Strategy~$3$. Figure~\ref{stratg3_u} gives a representation of the optimal 
control variables $u_1(t)$ and $u_2(t)$ of Strategy~$3$.
The minimum of the $J(u)$ functional \eqref{cost_func_strat3} 
of Strategy~3 is 0,94.
% -------------------------------
% -------------------------------
% -------------------------------
% -------------------------------
% -------------------------------

% -------------------------------

\section{Discussion}
\label{subsec:3.5}

\begin{figure}
\centering
\includegraphics[width=10cm]{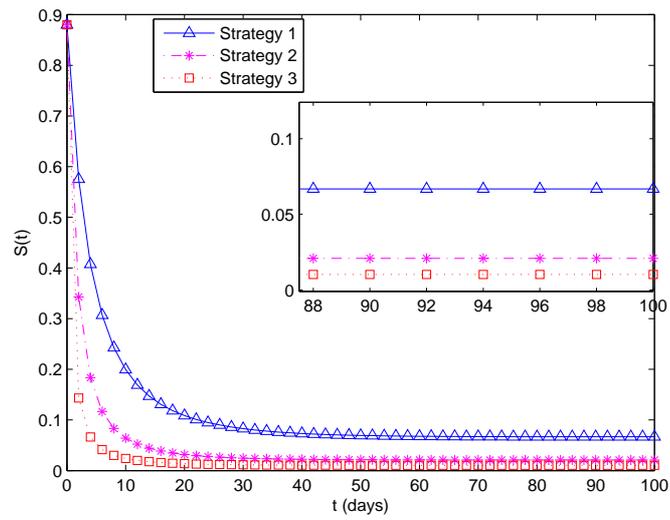}
\caption{Comparison between the curves of susceptible individuals $S(t)$
in case of Strategy~$1$, Strategy~$2$, and Strategy~$3$. \label{stratg123_S}}
\end{figure}

\begin{figure}
\centering
\includegraphics[width=10cm]{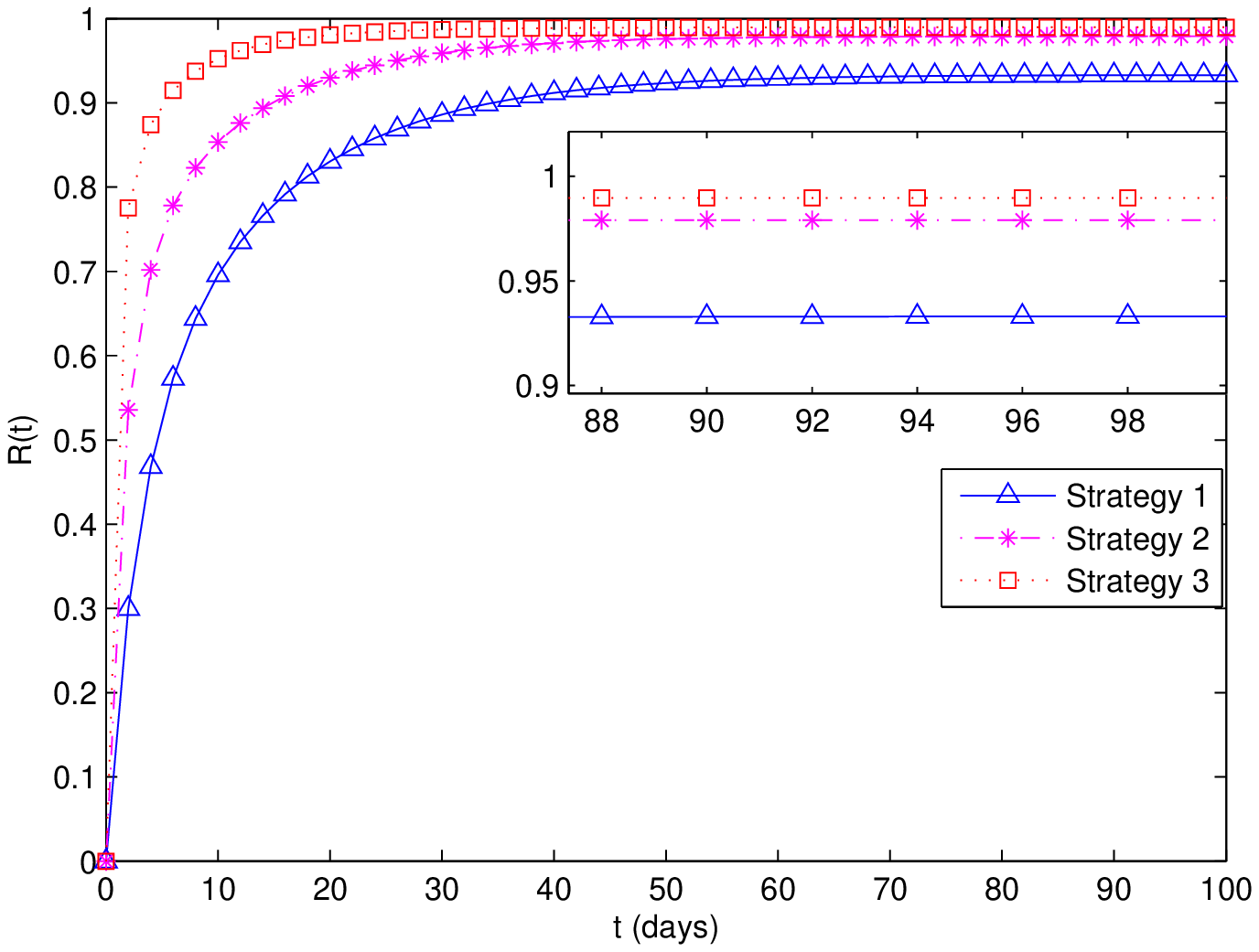}
\caption{Comparison between the curves of recovered individuals $R(t)$
in case of Strategy~$1$, Strategy~$2$, and Strategy~$3$. \label{stratg123_R}}
\end{figure}

\begin{figure}
\centering
\includegraphics[width=10cm]{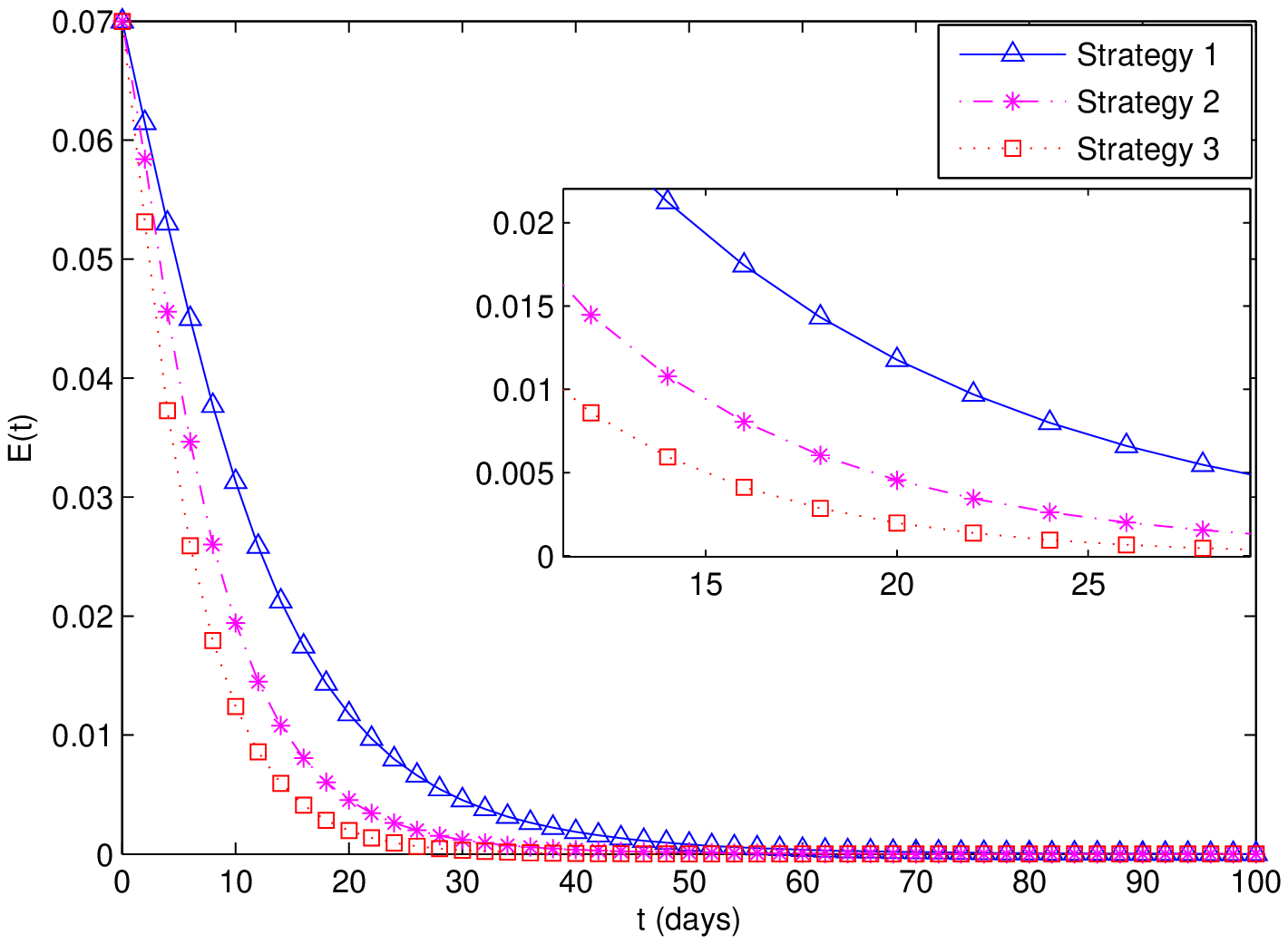}
\caption{Comparison between the curves of exposed individuals $E(t)$
in case of Strategy~$1$, Strategy~$2$, and Strategy~$3$. \label{stratg123_E}}
\end{figure}

\begin{figure}
\centering
\includegraphics[width=10cm]{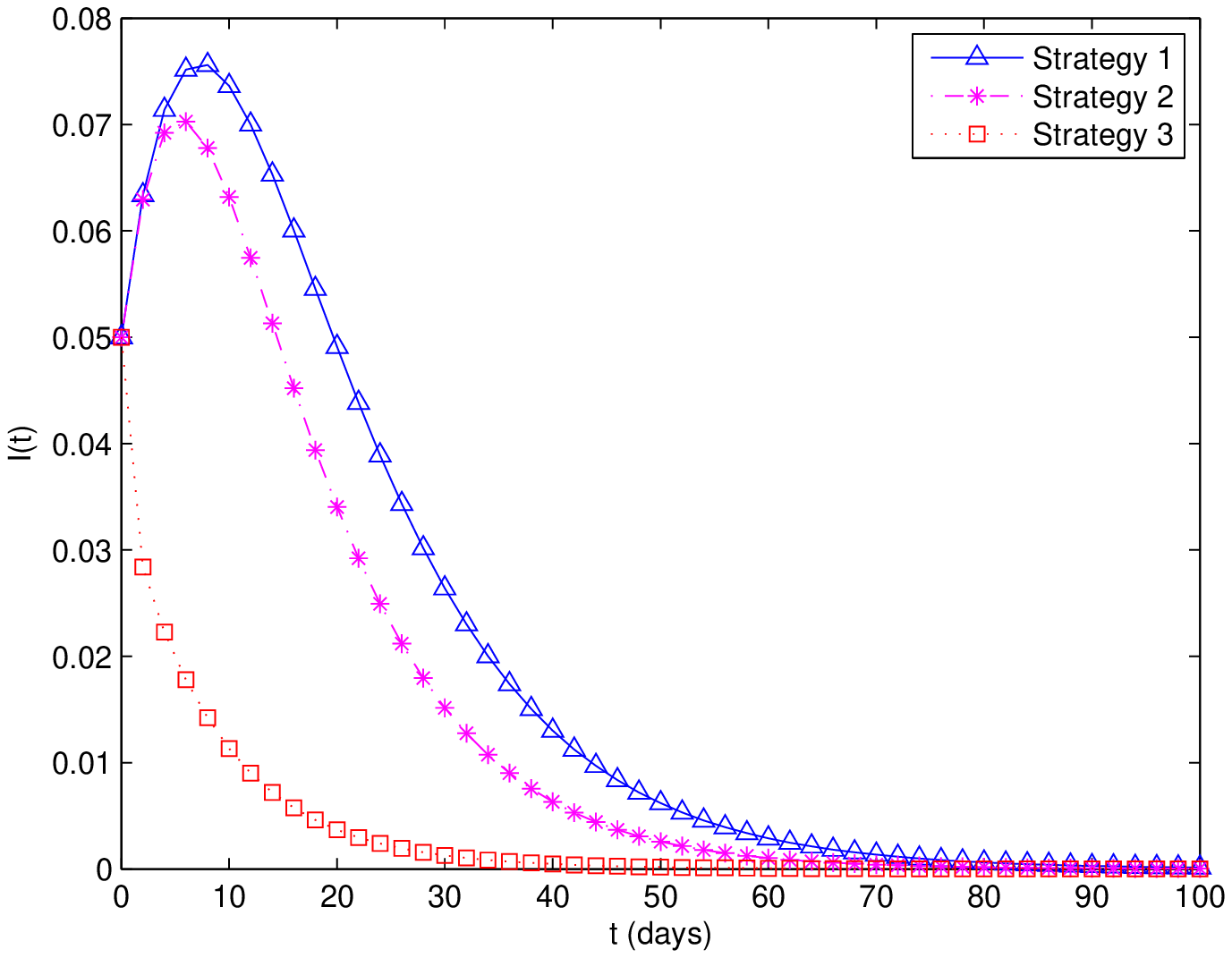}
\caption{Comparison between the curves of infected individuals $I(t)$
in case of Strategy~$1$, Strategy~$2$, and Strategy~$3$. \label{stratg123_I}}
\end{figure}

In this section we compare between the three strategies
of Section~\ref{sec:3}, and we discuss the obtained results.
Figures~\ref{stratg123_S} and \ref{stratg123_R} represent, respectively, the
time-dependent curve of susceptible $S$ and recovered individuals $R$.
The curves of $S$ and $R$ of Strategy~$3$ are slightly greater than the curves
of $S$ and $R$ of Strategy~$2$, which is logical because in the cost functional
of Strategy~$2$ we minimize the infected and the exposed individuals 
with vaccination, where in the cost functional of Strategy~$3$ we minimize the 
infected number with the cost of educational campaigns 
and vaccination. It is important in the comparison between the three curves to
note that the number of exposed and infected individuals decrease the most 
rapidly in case of Strategy~$3$ (see Figures~\ref{stratg123_E} 
and \ref{stratg123_I}). Moreover, there is not any peak in case of Strategy~$3$.
The results show the efficiency of an educational campaign in controlling Ebola 
virus when we couple it with the vaccination treatment control (as given in 
Strategy~$3$). We conclude that one can improve vaccination by educational 
campaigns, which has a great importance in poor countries who do not have 
the capacity to defend themselves against the virus.
% -------------------------------
% -------------------------------
% -------------------------------
% -------------------------------

% ------------------------------------

\section{Conclusion}
\label{sec:4}

We studied several optimal control problems 
for the SEIR model recently discussed by Rachah and Torres
in \cite{symcomp}, which provides a good description of the 2014
Ebola outbreak in West Africa. Precisely, we introduced 
a control $u(t)$ representing the vaccination 
rate at time $t$ (Section~\ref{subsec:3.1}). Then
we addressed the problem of reducing not only the number of
infected individuals $I(t)$ but also the number 
of exposed individuals $E(t)$ (Section~\ref{subsec:3.2}).
Finally, we investigated the integration of
an educational campaign about the virus into the population
(Section~\ref{subsec:3.3}).
It has been shown that an educational campaign has a great importance
with the vaccination treatment, mainly in countries that do not have
the capacity to defend themselves against the virus. As future work,
we plan to include in our study other factors. For instance,
we intend to include in the mathematical model a quarantine procedure.

% ------------------------------------

\begin{acknowledgements}
This research was partially supported by the
Institute of Mathematics of Toulouse, France (Rachah);
and by the Portuguese Foundation for Science and Technology (FCT),
within CIDMA project UID/MAT/04106/2013 and OCHERA project
PTDC/ EEI-AUT/1450/2012, co-financed by FEDER under POFC-QREN
with COMPETE reference FCOMP-01-0124-FEDER-028894 (Torres).
The authors would like to thank two referees for valuable
comments and helpful suggestions.
\end{acknowledgements}

% ------------------------------------------------------

% ------------------------------------

\end{document}